# EFFICIENT SCHEMES WITH UNCONDITIONALLY ENERGY STABILITY FOR THE ANISOTROPIC CAHN-HILLIARD EQUATION USING THE STABILIZED-SCALAR AUGMENTED VARIABLE (S-SAV) APPROACH

XIAOFENG YANG

**Abstract.** In this paper, we consider numerical approximations for the anisotropic Cahn-Hilliard equation. The main challenge of constructing numerical schemes with unconditional energy stabilities for this model is how to design proper temporal discretizations for the nonlinear terms with the strong anisotropy. We propose two, second order time marching schemes by combining the recently developed SAV approach with the linear stabilization approach, where three linear stabilization terms are added. These terms are shown to be crucial to remove the oscillations caused by the anisotropic coefficients, numerically. The novelty of the proposed schemes is that all nonlinear terms can be treated semi-explicitly, and one only needs to solve three decoupled linear equations with constant coefficients at each time step. We further prove the unconditional energy stabilities rigorously, and present various 2D and 3D numerical simulations to demonstrate the stability and accuracy.

## 1. Introduction

The classical Cahn-Hilliard equation (isotropic) is a typical equation of the diffusive phase field method. Its governing system is derived from the variational approach of the total free energy that usually includes two parts, a linear part (gradient entropy) and a nonlinear part (Ginzburg-Landau double well potential). This model has been extensively studied and applied to resolve the motion of free interfaces between multiple material components. About the most recent advances in modeling, algorithms or computational technologies, we refer to [1, 3, 4, 8, 11–15, 30, 31] and the references therein.

In this paper, we consider numerical approximations for solving the anisotropic Cahn-Hilliard equation, which was proposed in [21, 22], where a sufficiently big anisotropic coefficient is introduced to simulate the formation of faceted pyramids on nanoscale crystal surfaces. Comparing to the isotropic model, in addition to the stiffness issue from the interfacial width, a particular speciality of the anisotropic system is that the strong anisotropic coefficient $\gamma(\frac{\nabla \phi}{|\nabla \phi|})$ can induce large oscillations numerically, that increases the complexity for algorithm developments to a large extent. Therefore, there are very few successful attempts in designing efficient and energy stable schemes for this model. In [21, 22], the authors used the fully implicit method to discretize the nonlinear terms. Its energy stability is not provable, and the computational cost is high due to its nonlinear nature. In [2], a linear, first order scheme is developed based on the linear stabilization approach where all nonlinear terms are treated explicitly. Some extra linear stabilizers are added to enhance the stability. However, the energy stability is not provable for the anisotropic case even though the scheme is quite efficient and stable that allows for the large time step in computations. In [5], the

---

authors applied the convex splitting approach for solving this model, however, their scheme is only provably energy stable for the isotropic case. In [17], the authors proposed a second order scheme that is actually the modified version of scheme in [2], where the second order extrapolation is pre-estimated by the solutions of the first order scheme. However, such a stabilized predictor-corrector scheme is not provably energy stable as well.

Therefore, in this paper, the main purpose is to develop some efficient and effective numerical schemes to solve the anisotropic Cahn-Hilliard model. We expect that our schemes can combine the following three desired properties, i.e., (i) accurate (second order in time); (ii) stable (the unconditional energy dissipation law holds); and (iii) easy to implement and efficient (only need to solve some fully linear equations at each time step). To this end, we combine the the newly developed Scalar-Augmented-Variable (SAV) approach (cf. [18]), that is built upon the so-called Invariant-Energy-Quadratization (IEQ) approach (cf. [6,9,25–29,32–34]), with the linear stabilization approach (cf. [2,19]). The idea of SAV (or IEQ) method is to transform the total free energy integral (or integrand) into a quadratic function of a new, scalar auxiliary variable via a change of variables. Then, for the reformulated model in terms of the new variable that still retains the identical energy dissipation law, all nonlinear terms can be treated semi-explicitly. Due to the specialty of the new variable that is non-local type formally, one only needs to solve several decoupled, linear elliptic equations with constant coefficients at each time step. Meanwhile, since the anisotropic coefficient can cause large oscillations, we add some extra linear stabilizing terms. Numerical examples show that the combinations of these terms can efficiently remove oscillations and thus allow large time steps. We further show that the developed stabilized-SAV schemes are unconditional energy stable, and present various numerical examples to demonstrate the accuracy and stability. To the best of the author's knowledge, the schemes developed in this paper are the first such schemes for the anisotropic Cahn-Hilliard system that are theoretically proved to have all desired properties mentioned above.

The rest of the paper is organized as follows. In Section 2, we give a brief introduction of the governing PDE system for the anisotropic Cahn-Hilliard model. In Section 3, we develop two numerical schemes with second order accuracy for simulating the model, and rigorously prove their unconditional energy stabilities. Various 2D and 3D numerical experiments are given in Section 4 to validate the accuracy and efficiency of the proposed numerical schemes. Finally, some concluding remarks are given in Section 5.

## 2. Model equations and its energy law

Now we give a brief description for the anisotropic Cahn-Hilliard equation [21, 22]. Let $\Omega$ be a smooth, open, bounded, connected domain in $\mathbb{R}^d, d = 2, 3$. Let $\phi$ be an order parameter which takes the values $\pm 1$ in the two phases with a smooth transitional layer of thickness $\epsilon$. We consider the total free energy as follows,

$$(2.1) \qquad \mathcal{E}(\phi) = \int_\Omega \gamma(\boldsymbol{n})\big(\frac{1}{2}|\nabla\phi|^2 + \frac{1}{\epsilon^2}F(\phi)\big) + \frac{\beta}{2}G(\phi)d\boldsymbol{x},$$

where $\gamma(\boldsymbol{n})$ is a function describing the anisotropic property and $\boldsymbol{n}$ is the interface normal defined as $\boldsymbol{n} = \frac{\nabla\phi}{|\nabla\phi|}$. Therefore, for 2D,

$$(2.2) \qquad \boldsymbol{n} = (n_1, n_2)^T = \frac{1}{\sqrt{\phi_x^2 + \phi_y^2}}(\phi_x, \phi_y)^T,$$

and for 3D,

$$(2.3) \quad \boldsymbol{n} = (n_1, n_2, n_3)^T = \frac{1}{\sqrt{\phi_x^2 + \phi_y^2 + \phi_z^2}}(\phi_x, \phi_y, \phi_z)^T.$$

The nonlinear potential $F(\phi)$ takes the usual Ginzburg-Landau double well potential that reads as

$$(2.4) \quad F(\phi) = \frac{1}{4}(\phi^2 - 1)^2.$$

The anisotropic function may takes the fourfold form that reads as

$$(2.5) \quad \gamma(\boldsymbol{n}) = 1 + \alpha\cos(4\theta) = 1 + \alpha(4\sum_{i=1}^{d} n_i^4 - 3),$$

where $\theta$ denotes the orientation angle of the interfacial normal to the interface. The non-negative parameter $\alpha$ in (2.5) describes the intensity of anisotropy. When $\alpha = 0$, the system degenerates to the isotropic Cahn-Hilliard model. It is indicated in [21, 22], that a sufficient big $\alpha$ would produce a strongly anisotropic system, i.e., the underline Cahn-Hilliard equation is ill-posed. In order to regularize the system, an extra potential $G(\phi)$ is added to penalize infinite curvatures in the resulting corners and $\beta$ is the magnitude of the regularization parameter. Two kinds of regularization terms are generally considered. The first is the linear regularization based on the bi-Laplacian of the phase variable, that reads as

$$(2.6) \quad G_{linear}(\phi) = (\Delta\phi)^2,$$

and the second is the nonlinear Willmore regularization that reads as

$$(2.7) \quad G_{will}(\phi) = (\Delta\phi - \frac{1}{\epsilon^2}f(\phi))^2.$$

where $f(\phi) = F'(\phi) = \phi(\phi^2 - 1)$.

By taking the $H^{-1}$ gradient flow on the total free energy $\mathcal{E}(\phi)$, we obtain the anisotropic Cahn-Hilliard system with the linear regularization:

$$(2.8) \quad \phi_t = \nabla \cdot (M(\phi)\nabla\mu),$$

$$(2.9) \quad \mu = \frac{1}{\epsilon^2}\gamma(\boldsymbol{n})f(\phi) - \nabla \cdot \mathbf{m} + \beta\Delta^2\phi;$$

and the strongly anisotropic Cahn-Hilliard system with the Willmore regularization:

$$(2.10) \quad \phi_t = \nabla \cdot (M(\phi)\nabla\mu),$$

$$(2.11) \quad \mu = \frac{1}{\epsilon^2}\gamma(\boldsymbol{n})f(\phi) - \nabla \cdot \mathbf{m} + \beta(\Delta - \frac{1}{\epsilon^2}f'(\phi))(\Delta\phi - \frac{1}{\epsilon^2}f(\phi)),$$

where $f'(\phi) = 3\phi^2 - 1$, $M(\phi) \geq 0$ is the mobility function that depends on the phase variable $\phi$. The vector field $\mathbf{m}$ is defined as

$$(2.12) \quad \mathbf{m} = \gamma(\boldsymbol{n})\nabla\phi + \frac{\mathbb{P}\nabla_{\boldsymbol{n}}\gamma(\boldsymbol{n})}{|\nabla\phi|}(\frac{1}{2}|\nabla\phi|^2 + \frac{1}{\epsilon^2}F(\phi)),$$

where the matrix $\mathbb{P} = \mathbb{I} - \boldsymbol{n}\boldsymbol{n}^T$.

Without the loss of generality, we adopt the periodic boundary condition to remove all complexities associated with the boundary integrals in this study. We remark that the boundary conditions can also be the no-flux type as

$$\frac{\partial \phi}{\partial \mathbf{n}}\Big|_{\partial\Omega} = \frac{\partial \mu}{\partial \mathbf{n}}\Big|_{\partial\Omega} = \frac{\partial \omega}{\partial \mathbf{n}}\Big|_{\partial\Omega} = 0, \tag{2.13}$$

where $w = \Delta\phi$ for linear regularization model and $w = \Delta\phi - \frac{1}{\epsilon^2}f(\phi)$ for the Willmore regularization model, $\mathbf{n}$ is the outward normal of the computational domain $\Omega$. All numerical analysis in this paper can be carried out to the no-flux boundary conditions without any further difficulties.

The model equations (2.8)-(2.9) and (2.10)-(2.11) follow the dissipative energy law. By taking the $L^2$ inner product of (2.8) with $-\mu$, and of (2.9) with $\phi_t$, using the integration by parts and combining the obtained two equalities, we obtain we obtain

$$\frac{d}{dt}\mathcal{E}(\phi) = -\|\sqrt{M(\phi)}\nabla\mu\|^2 \leq 0. \tag{2.14}$$

Meanwhile, the Cahn-Hilliard model conserves the local mass density. By taking the $L^2$ inner product of (2.8) with 1, one can obtain the mass conservation property directly using integration by parts, that reads as

$$\frac{d}{dt}\int_\Omega \phi d\boldsymbol{x} = 0. \tag{2.15}$$

## 3. Numerical Schemes

We develop in this section a set of linear, second-order, unconditionally energy stable schemes for solving the anisotropic model (2.8)-(2.9) and (2.10)-(2.11). The main challenging issues are to develop suitable approaches to discretize the nonlinear terms, particularly, the terms associated with the anisotropic coefficient $\gamma(\boldsymbol{n})$, as well as the nonlinear term from the Willmore regularization potential (the term associated with $\beta$ in (2.11)).

For the simple phase field models, e.g., the isotropic Cahn-Hilliard equation or Allen-Cahn equation where the only numerical challenge is to discretize the cubic term from the double-well potential, there are many techniques to preserve the unconditional energy stability. Among them, two prevalent methods have been studied extensively. The first one is the so-called convex splitting approach [7, 16], where the convex part of the potential is treated implicitly and the concave part is treated explicitly. The convex splitting approach is energy stable, however, it usually produces a nonlinear scheme at most cases, thus the implementation is complicated and the computational cost is high. The second one is the so-called stabilization approach [10, 17, 19, 20, 23, 24], where the nonlinear term is treated in the simple explicit way. In order to preserve the energy stability, a linear stabilizing term has to be added, and the magnitude of that term usually depends on the upper bound of the second order derivative of the double-well potential. The stabilization approach introduces a purely linear scheme, thus it is absolutely easy to implement. However, if there does not exist any finite upper bound for that second order derivative, one must make some reasonable revisions to the nonlinear potential in order to obtain a finite bound, for example, a quadratic cut-off function for the double-well potential. Such a method is particularly reliable for those models satisfying the maximum principle. Otherwise, if the maximum principle does not hold, the revisions to the nonlinear potentials may lead to spurious solutions.

It is remarkable that the ideas of the convex splitting or the stabilization approach are not suitable for the anisotropic model. The stabilization approach cannot provide ideal schemes, since (i) it is uncertain whether the PDE solution could satisfy any maximum principle; (ii) it is extremely hard obtain the energy stability for the developed scheme if all nonlinear terms are treated explicitly.

Meanwhile, it is particularly not clear and questionable whether the nonlinear potentials multiplied with the anisotropic coefficient could be split into the combinations of the convex and concave parts, which exclude the convex splitting approach.

We aim to develop some more effective and efficient numerical schemes in this paper. More precisely, we expect that the schemes are efficient (linear system), stable (unconditionally energy stable), and accurate (second order). To this end, we developed the so-called stabilized-SAV approach by combining the recently SAV approach with the linear stabilization approach mentioned above, without worrying about whether the continuous/discrete maximum principle holds or a convexity/concavity splitting exists. First, we transform the total free energy integral into a quadratic function of a new, scalar auxiliary variable. For the reformulated model, we treat all nonlinear terms in the semi-explicit way. Furthermore, some crucial linear stabilizing terms are added. These terms cannot only enhance the stabilities, but also are particularly efficient to remove oscillations caused by the anisotropic coefficient.

### 3.1. Linear regularization model.

We first deal with the linear regularization model. Let us define an auxiliary variable as follows:

$$\text{(3.1)} \quad U = \sqrt{\int_\Omega \gamma(\boldsymbol{n})\big(\tfrac{1}{2}|\nabla\phi|^2 + \tfrac{1}{\epsilon^2}F(\phi)\big)d\boldsymbol{x} + B},$$

where $B$ is any constant that ensures the radicand positive (in all numerical examples, we let $B = 1$). Thus the total free energy (2.1) can be rewritten as

$$\text{(3.2)} \quad \mathscr{E}(U, \phi) = U^2 - B + \frac{\beta}{2}\int_\Omega \Delta\phi^2 d\boldsymbol{x},$$

Using these new variables $U$, we then have the following equivalent PDE system,

$$\text{(3.3)} \quad \phi_t = \nabla \cdot (M(\phi)\nabla\mu),$$
$$\text{(3.4)} \quad \mu = H(\phi)U + \beta\Delta^2\phi.$$

By taking the time derivative for the new variable $U$, we obtain

$$\text{(3.5)} \quad U_t = \frac{1}{2}\int_\Omega H(\phi)\phi_t d\boldsymbol{x},$$

where

$$\text{(3.6)} \quad H(\phi) = \frac{\tfrac{1}{\epsilon^2}\gamma(\boldsymbol{n})f(\phi) - \nabla \cdot \boldsymbol{m}}{\sqrt{\int_\Omega \gamma(\boldsymbol{n})\big(\tfrac{1}{2}|\nabla\phi|^2 + \tfrac{1}{\epsilon^2}F(\phi)\big)d\boldsymbol{x} + B}}.$$

The transformed system (3.3)-(3.5) forms a closed PDE system with the following initial conditions,

$$\text{(3.7)} \quad \begin{cases} \phi(t=0) = \phi^0, \\ U(t=0) = \sqrt{\int_\Omega \gamma(\boldsymbol{n}^0)\big(\tfrac{1}{2}|\nabla\phi^0|^2 + \tfrac{1}{\epsilon^2}F(\phi^0)\big)d\boldsymbol{x} + B}. \end{cases}$$

The system (3.3)-(3.5) also follows an energy dissipative laws in terms of the new variable $U$ and $\phi$. By taking the $L^2$ inner product of (3.3) with $-\mu$, of (3.4) with $\phi_t$, multiplying (3.5) with $2U$, performing integration by parts and summing all equalities up, we can obtain the energy dissipation law of the new system (3.3)-(3.5) as

$$\frac{d}{dt}\mathscr{E}(U,\phi) = -\|\sqrt{M(\phi)}\nabla\mu\|^2 \leq 0. \tag{3.8}$$

We note that the new transformed system (3.3)-(3.5) is equivalent to the original system (2.8)-(2.9) for the time continuous case since (3.1) can be easily obtained by integrating (3.5) with respect to the time. Next we will develop unconditionally energy stable linear numerical schemes for time stepping of the transformed system (3.3)-(3.5), and the proposed schemes should formally follow the new energy dissipation law (3.8) in the discrete sense, instead of the energy law for the originated system (2.14). Even though it is shown by many numerical examples in section 4 that these two free energies, (3.2) and (2.1), are quite different quantitatively after long time computations, the developed schemes that follow the new energy dissipation law (3.8) can still ensure the original energy decayed, numerically.

Let $\delta t > 0$ be a time step size and set $t^n = n\delta t$ for $0 \leq n \leq N$ with $T = N\delta t$. We also denote the $L^2$ inner product of any two spatial functions $f_1(\boldsymbol{x})$ and $f_2(\boldsymbol{x})$ by $(f_1(\boldsymbol{x}), f_2(\boldsymbol{x})) = \int_\Omega f_1(\boldsymbol{x})f_2(\boldsymbol{x})d\boldsymbol{x}$, and the $L^2$ norm of the function $f(\boldsymbol{x})$ by $\|f\|^2 = (f,f)$. Let $\psi^n$ denotes the numerical approximation to $\psi(\cdot,t)|_{t=t^n}$ for any function $\psi$.

We construct a second order numerical scheme based on the second order backward differentiation formula (BDF2).

**Scheme 1.** *Assuming $\phi^n, U^n$ and $\phi^{n-1}, U^{n-1}$ are known, we update $\phi^{n+1}, U^{n+1}$ by solving*

$$\frac{3\phi^{n+1} - 4\phi^n + \phi^{n-1}}{2\delta t} = \nabla(M(\phi^{\star,n+1})\nabla\mu^{n+1}), \tag{3.9}$$

$$\mu^{n+1} = H^{\star,n+1}U^{n+1} + \beta\Delta^2\phi^{n+1} \tag{3.10}$$
$$+ \frac{S_1}{\epsilon^2}(\phi^{n+1} - 2\phi^n + \phi^{n-1}) - S_2\Delta(\phi^{n+1} - 2\phi^n + \phi^{n-1})$$
$$+ S_3\Delta^2(\phi^{n+1} - 2\phi^n + \phi^{n-1}),$$

$$3U^{n+1} - 4U^n + U^{n-1} = \frac{1}{2}\int_\Omega H^{\star,n+1}(3\phi^{n+1} - 4\phi^n + \phi^{n-1})d\boldsymbol{x}, \tag{3.11}$$

*where*

$$\phi^{\star,n+1} = 2\phi^n - \phi^{n-1}, \ H^{\star,n+1} = H(\phi^{\star,n+1}), \tag{3.12}$$

*and $S_{i,i=1,2,3}$ are three positive stabilizing parameters.*

**Remark 3.1.** *Three second order linear stabilizers (associated with $S_1, S_2, S_3$) are added in the scheme. The first two stabilizers, $\frac{S_1}{\epsilon^2}\Delta(\phi^{n+1} - 2\phi^n + \phi^{n-1})$ and $-S_2\Delta(\phi^{n+1} - 2\phi^n + \phi^{n-1})$ are two commonly used linear stabilizers in the linear stabilization method for solving the isotropic or anisotropic phase field model, (cf. [2] for the anisotropic model, and [19] for the isotropic model). The errors that these two terms introduce are of order $\frac{S_1}{\epsilon^2}\delta t^2\phi_{tt}(\cdot)$ and $S_2\delta t^2\Delta\phi_{tt}(\cdot)$, respectively, which are of the same order as the error introduced by the second order extrapolation of the nonlinear term $f(\phi)$ and the Laplacian term $\Delta\phi$. The last stabilizing term associated with $S_3$ is of the order $S_3\delta t^2\Delta^2\phi_{tt}(\cdot)$. In practice, we choose $S_3 = O(\beta)$. Numerical examples show that the combinations of the three stabilizers are crucial to remove all oscillations induced by the anisotropic coefficient $\gamma(\boldsymbol{n})$, since the term $\boldsymbol{n} = \frac{\nabla\phi}{|\nabla\phi|}$ changes its sign frequently as long as when $|\nabla\phi|$ is close to zero (cf. Fig. 4.4 in section 4).*

Apparently, one needs to solve a nonlocal, coupled system for $\phi^{n+1}$ and $U^{n+1}$ in the scheme (3.9)-(3.11). But actually, in practice, if the mobility $M(\phi)$ is a constant, e.g., $M(\phi) = M_0$, we can simplify the scheme through the following procedure.

We first rewrite (3.11) as follows,

$$(3.13) \qquad U^{n+1} = \frac{1}{2}\int_\Omega H^{\star,n+1}\phi^{n+1} + g^n,$$

where $g^n = \frac{4U^n - U^{n-1}}{3} - \frac{1}{2}\int_\Omega H^\star \frac{4\phi^n - \phi^{n-1}}{3}d\boldsymbol{x}$. Then the scheme (3.9) can be written as

$$(3.14) \quad (\frac{3}{2M_0\delta t} - \frac{S_1}{\epsilon^2}\Delta + S_2\Delta^2 - (S_3+\beta)\Delta^3)\phi^{n+1} - \frac{1}{2}\Delta H^{\star,n+1}\int_\Omega H^{\star,n+1}\phi^{n+1}d\boldsymbol{x} = \tilde{g}^n,$$

where

$$(3.15) \qquad \tilde{g}^n = \frac{4\phi^n - \phi^{n-1}}{2M_0\delta t} + \Delta(H^{\star,n+1}g^n - \frac{S_1}{\epsilon^2}\phi^{\star,n+1} + S_2\Delta\phi^{\star,n+1} - S_3\Delta^2\phi^{\star,n+1}).$$

Define an linear operator $\chi^{-1}(\cdot)$, such that for any $\phi \in L^2(\Omega)$, $\psi = \chi^{-1}(\phi)$ is defined as

$$(3.16) \qquad (\frac{3}{2M_0\delta t} - \frac{S_1}{\epsilon^2}\Delta + S_2\Delta^2 - (S_3+\beta)\Delta^3)\psi = \phi.$$

By applying the operator $\chi^{-1}$ to (3.14), then we obtain

$$(3.17) \qquad \phi^{n+1} - \frac{1}{2}\chi^{-1}(\Delta H^{\star,n+1})\int_\Omega H^{\star,n+1}\phi^{n+1}d\boldsymbol{x} = \chi^{-1}(\tilde{g}^n).$$

By taking the $L^2$ inner product with $H^{\star,n+1}$, we obtain

$$(3.18) \qquad \int_\Omega H^{\star,n+1}\phi^{n+1}d\boldsymbol{x} = \frac{\int_\Omega H^{\star,n+1}\chi^{-1}(\tilde{g}^n)d\boldsymbol{x}}{1 - \frac{1}{2}\int_\Omega H^{\star,n+1}\chi^{-1}(\Delta H^{\star,n+1})d\boldsymbol{x}}.$$

Note the term in the denominator $-\int_\Omega H^{\star,n+1}\chi^{-1}(\Delta H^{\star,n+1})d\boldsymbol{x} \geq 0$ since $-\chi^{-1}\Delta$ is a positive definite operator. Therefore, in the computations, one only needs to find $\psi_1 = \chi^{-1}(\tilde{g}^n)$ and $\psi_2 = \chi^{-1}(\Delta H^{\star,n+1})$, that means to solve the following two sixth order sub-equations,

$$(3.19) \qquad (\frac{3}{2M_0\delta t} - \frac{S_1}{\epsilon^2}\Delta + S_2\Delta^2 - (S_3+\beta)\Delta^3)\psi_1 = f^n,$$

and

$$(3.20) \qquad (\frac{3}{2M_0\delta t} - \frac{S_1}{\epsilon^2}\Delta + S_2\Delta^2 - (S_3+\beta)\Delta^3)\psi_2 = \Delta H^{\star,n+1},$$

with the periodic boundary conditions. Once $\psi_1$ and $\psi_2$ are obtained, by applying (3.18) to get $\int_\Omega H^{\star,n+1}\phi^{n+1}d\boldsymbol{x}$, we then solve the third sixth order equation (3.14) to obtain $\phi^{n+1}$.

To summarize, the scheme (3.9)-(3.11) can be easily implemented in the following manner:

- Compute $\psi_1$ and $\psi_2$ by solving two sixth-order equations with constant coefficients, (3.19) and (3.20);
- Compute $\int_\Omega H^{\star,n+1}\phi^{n+1}$ from (3.18) and $U^{n+1}$ from (3.13);
- Compute $\phi^{n+1}$ by solving the third sixth-order equation with constant coefficients (3.14).

Hence, the total cost at each time step are essentially solving three sixth-order equations with constant coefficients. We note that these sixth-order equations with periodic boundary conditions can be easily computed, hence, this scheme is extremely efficient and easy to implement.

Now we prove the scheme (3.9)-(3.11) is unconditionally energy stable as follows.

**Theorem 3.1.** *The scheme (3.9)-(3.11) is unconditionally energy stable which satisfies the following discrete energy dissipation law,*

$$\frac{1}{\delta t}(E_{linear}^{n+1} - E_{linear}^n) \leq -\|\sqrt{M(\phi^{\star,n+1})}\nabla\mu^{n+1}\|^2 \leq 0, \tag{3.21}$$

*where*

$$\begin{aligned}
E_{linear}^{n+1} =& \frac{(U^{n+1})^2 + (2U^{n+1} - U^n)^2}{2} + \frac{\beta}{2}\Big(\frac{\|\Delta\phi^{n+1}\|^2 + \|2\Delta\phi^{n+1} - \Delta\phi^n\|^2}{2}\Big) \\
& + \frac{S_1}{\epsilon^2}\frac{\|\phi^{n+1} - \phi^n\|^2}{2} + S_2\frac{\|\nabla\phi^{n+1} - \nabla\phi^n\|^2}{2} + S_3\frac{\|\Delta\phi^{n+1} - \Delta\phi^n\|^2}{2}.
\end{aligned} \tag{3.22}$$

*Proof.* By taking the $L^2$ inner product of (3.9) with $-2\delta t\mu^{n+1}$, we obtain

$$-(3\phi^{n+1} - 4\phi^n + \phi^{n-1}, \mu^{n+1}) = 2\delta t\|\sqrt{M(\phi^{\star,n+1})}\nabla\mu^{n+1}\|^2. \tag{3.23}$$

By taking the $L^2$ inner product of (3.10) with $3\phi^{n+1} - 4\phi^n + \phi^{n-1}$, and using integration by parts, we obtain

$$\begin{aligned}
(\mu^{n+1}, 3\phi^{n+1} - 4\phi^n + \phi^{n-1}) =& U^{n+1}(H^{\star,n+1}, 3\phi^{n+1} - 4\phi^n + \phi^{n-1}) \\
& + \frac{S_1}{\epsilon^2}(\phi^{n+1} - 2\phi^n + \phi^{n-1}, 3\phi^{n+1} - 4\phi^n + \phi^{n-1}) \\
& + S_2(\nabla(\phi^{n+1} - 2\phi^n + \phi^{n-1}), \nabla(3\phi^{n+1} - 4\phi^n + \phi^{n-1})) \\
& + S_3(\Delta(\phi^{n+1} - 2\phi^n + \phi^{n-1}), \Delta(3\phi^{n+1} - 4\phi^n + \phi^{n-1})) \\
& + \beta(\Delta\phi^{n+1}, \Delta(3\phi^{n+1} - 4\phi^n + \phi^{n-1})).
\end{aligned} \tag{3.24}$$

By multiplying (3.11) with $-2U^{n+1}$, we obtain

$$-2(3U^{n+1} - 4U^n + U^{n-1})U^{n+1} = -U^{n+1}\int_\Omega H^{\star,n+1}(3\phi^{n+1} - 4\phi^n + \phi^{n-1})d\boldsymbol{x} \tag{3.25}$$

Combining the above equations and applying the following two identities

$$\begin{aligned}
2a(3a - 4b + c) &= a^2 + (2a - b)^2 - b^2 - (2b - c)^2 + (a - 2b + c)^2, \\
(3a - 4b + c)(a - 2b + c) &= (a - b)^2 - (b - c)^2 + 2(a - 2b + c)^2,
\end{aligned} \tag{3.26}$$

we obtain

$$\begin{aligned}
& \Big((U^{n+1})^2 + (2U^{n+1} - U^n)^2\Big) - \Big((U^{n+1})^2 + (2U^n - U^{n-1})^2\Big) \\
& + \frac{\beta}{2}\Big(\|\Delta\phi^{n+1}\|^2 + \|2\Delta\phi^{n+1} - \Delta\phi^n\|^2\Big) - \frac{\beta}{2}\Big(\|\Delta\phi^n\|^2 + \|2\Delta\phi^n - \Delta\phi^{n-1}\|^2\Big) \\
& + \frac{S_1}{\epsilon^2}\|\phi^{n+1} - \phi^n\|^2 - s\|\phi^n - \phi^{n-1}\|^2 \\
& + S_2\|\nabla\phi^{n+1} - \nabla\phi^n\|^2 - s\|\nabla\phi^n - \nabla\phi^{n-1}\|^2 \\
& + S_3\|\Delta\phi^{n+1} - \Delta\phi^n\|^2 - s\|\Delta\phi^n - \Delta\phi^{n-1}\|^2 \\
& + (U^{n+1} - 2U^n + U^{n-1})^2 + \frac{2S_1}{\epsilon^2}\|\phi^{n+1} - 2\phi^n + \phi^{n-1}\|^2 \\
& + 2S_2\|\nabla(\phi^{n+1} - 2\phi^n + \phi^{n-1})\|^2 + 2S_3\|\Delta(\phi^{n+1} - 2\phi^n + \phi^{n-1})\|^2 + \frac{\beta}{2}\|\Delta(\phi^{n+1} - 2\phi^n + \phi^{n-1})\|^2 \\
=& -2\delta t\|\sqrt{M(\phi^{\star,n+1})}\nabla\mu^{n+1}\|^2.
\end{aligned}$$

Finally, we obtain the desired result after dropping some positive terms. □

**Remark 3.2.** *Heuristically, $\frac{1}{\delta t}(E_{linear}^{n+1} - E_{linear}^n)$ is a second-order approximation of $\frac{d}{dt}\mathcal{E}(\phi, U)$ at $t = t^{n+1}$. For any smooth variable $\psi$ with time, we have*

$$\frac{\|\psi^{n+1}\|^2 - \|2\psi^{n+1} - \psi^n\|^2}{2\delta t} - \frac{\|\psi^n\|^2 - \|2\psi^n - \psi^{n-1}\|^2}{2\delta t}$$

(3.27)
$$\cong \frac{\|\psi^{n+2}\|^2 - \|\psi^n\|^2}{2\delta t} + O(\delta t^2) \cong \frac{d}{dt}\|\psi(t^{n+1})\|^2 + O(\delta t^2),$$

*and*

(3.28)
$$\frac{\|\psi^{n+1} - \psi^n\|^2 - \|\psi^n - \psi^{n-1}\|^2}{2\delta t} \cong O(\delta t^2).$$

**Remark 3.3.** *It is also straightforward to develop the second order Crank-Nicolson scheme where the linear stabilizers terms still form like $\psi^{n+1} - 2\psi^n + \psi^{n-1}$. We omit the details to the interested readers since the proof of energy stability is quite similar to Theorem 3.1. In addition, although we consider only time discrete schemes in this study, the results can be carried over to any consistent finite-dimensional Galerkin approximations in the space since the proofs are all based on a variational formulation with all test functions in the same space as the space of the trial functions.*

### 3.2. Willmore regularization model.

We deal with the Willmore regularization model in this subsection. Similar to the linear regularization case, we define an auxiliary variable as follows:

(3.29)
$$U = \sqrt{\int_\Omega \Big(\gamma(\boldsymbol{n})\big(\frac{1}{2}|\nabla\phi|^2 + \frac{1}{\epsilon^2}F(\phi)\big) + \frac{\beta}{2}(\Delta\phi - \frac{1}{\epsilon^2}f(\phi))^2\Big)d\boldsymbol{x} + B},$$

where $B$ is any constant that ensures the radicand positive, thus the total free energy (2.1) can be rewritten as

(3.30)
$$\mathscr{E}(U, \phi) = U^2 - B$$

Using these new variables $U$, we then have the following equivalent PDE system,

(3.31) $$\phi_t = \nabla \cdot (M(\phi)\nabla\mu),$$

(3.32) $$\mu = Z(\phi)U,$$

(3.33) $$U_t = \frac{1}{2}\int_\Omega Z(\phi)\phi_t d\boldsymbol{x}.$$

where

(3.34)
$$Z(\phi) = \frac{\frac{1}{\epsilon^2}\gamma(\boldsymbol{n})f(\phi) - \nabla \cdot \boldsymbol{m} + \beta(\Delta - \frac{1}{\epsilon^2}f'(\phi))(\Delta\phi - \frac{1}{\epsilon^2}f(\phi))}{\sqrt{\int_\Omega \Big(\gamma(\boldsymbol{n})\big(\frac{1}{2}|\nabla\phi|^2 + \frac{1}{\epsilon^2}F(\phi)\big) + \frac{\beta}{2}(\Delta\phi - \frac{1}{\epsilon^2}f(\phi))^2\Big)d\boldsymbol{x} + B}}.$$

The system (3.31)-(3.33) is equipped with the following initial conditions,

(3.35)
$$\begin{cases} \phi(t=0) = \phi^0, \\ U(t=0) = \sqrt{\int_\Omega \Big(\gamma(\boldsymbol{n}^0)\big(\frac{1}{2}|\nabla\phi^0|^2 + \frac{1}{\epsilon^2}F(\phi^0)\big) + \frac{\beta}{2}(\Delta\phi^0 - \frac{1}{\epsilon^2}f(\phi^0))^2\Big)d\boldsymbol{x} + B}. \end{cases}$$

It is easy to see that the system (3.31)-(3.33) follows an energy dissipative laws in terms of the new variable $U$ and $\phi$. By taking the $L^2$ inner product of (3.3) with $-\mu$, of (3.4) with $\phi_t$, multiplying (3.5) with $2U$, performing integration by parts and summing all equalities up, we can obtain the energy dissipation law of the new system (3.3)-(3.5) as

$$\frac{d}{dt}\mathscr{E}(U,\phi) = -\|\sqrt{M(\phi)}\nabla\mu\|^2 \leq 0. \tag{3.36}$$

We present the second order numerical scheme based on the BDF2 for solving the transformed model (3.31)-(3.33).

**Scheme 2.** *Assuming $\phi^n, U^n$ and $\phi^{n-1}, U^{n-1}$ are known, we update $\phi^{n+1}, U^{n+1}$ by solving*

$$\frac{3\phi^{n+1} - 4\phi^n + \phi^{n-1}}{2\delta t} = \nabla(M(\phi^{\star,n+1})\nabla\mu^{n+1}), \tag{3.37}$$

$$\mu^{n+1} = Z^{\star,n+1}U^{n+1} + \frac{S_1}{\epsilon^2}(\phi^{n+1} - 2\phi^n + \phi^{n-1}) \tag{3.38}$$
$$- S_2\Delta(\phi^{n+1} - 2\phi^n + \phi^{n-1}) + S_3\Delta^2(\phi^{n+1} - 2\phi^n + \phi^{n-1}),$$

$$3U^{n+1} - 4U^n + U^{n-1} = \frac{1}{2}\int_\Omega Z^{\star,n+1}(3\phi^{n+1} - 4\phi^n + \phi^{n-1})d\boldsymbol{x}, \tag{3.39}$$

*where $Z^{\star,n+1} = Z(\phi^{\star,n+1})$, $S_{i,i=1,2,3}$ are three positive stabilizing parameters.*

Since the scheme (3.37)-(3.39) is almost identical to the scheme (3.9)-(3.11) formally, thus we omit the details of its implementations and only present the theorem of energy stability as follows.

**Theorem 3.2.** *The scheme (3.37)-(3.39) is unconditionally energy stable which satisfies the following discrete energy dissipation law,*

$$\frac{1}{\delta t}(E^{n+1}_{will} - E^n_{will}) \leq -\|\sqrt{M(\phi^{\star,n+1})}\nabla\mu^{n+1}\|^2 \leq 0, \tag{3.40}$$

*where*

$$E^{n+1}_{will} = \frac{(U^{n+1})^2 + (2U^{n+1} - U^n)^2}{2} + \frac{S_1}{\epsilon^2}\frac{\|\phi^{n+1} - \phi^n\|^2}{2} \tag{3.41}$$
$$+ S_2\frac{\|\nabla\phi^{n+1} - \nabla\phi^n\|^2}{2} + S_3\frac{\|\Delta\phi^{n+1} - \Delta\phi^n\|^2}{2}.$$

## 4. Numerical simulation

We now present various numerical examples to validate the proposed schemes and demonstrate their accuracy, energy stability and efficiency (particularly for removing oscillations). Here, we choose the periodic boundary conditions and set the computational domain as $\Omega = [0, 2\pi]^d, d = 2, 3$. We use the Fourier-spectral method to discretize the space, where $129 \times 129$ Fourier modes are used for 2D simulations, and $129 \times 129 \times 129$ Fourier modes are used for 3D simulations.

If not explicitly specified, the default values of order parameters and stabilization parameters are set as follows,

$$\epsilon = 6\mathrm{e}{-2}, \alpha = 0.3, \beta = 6\mathrm{e}{-4}, S_1 = 2, S_2 = 2, S_3 = 1\mathrm{e}{-3}. \tag{4.1}$$

**4.1. Accuracy test.** We first perform numerical simulations to test the convergence rates of the two proposed schemes, the scheme (3.9)-(3.11) for the linear regularization model (denoted by BDF2-L) and the scheme (3.37)-(3.39) for the Willmore regularization model (denoted by BDF2-W).

| $\delta t$ | BDF2-L | Order | BDF2-W | Order |
|---|---|---|---|---|
| 1e−2 | 4.01e−5 | | 3.93e−5 | |
| 5e−3 | 1.15e−5 | 1.80 | 1.01e−5 | 1.96 |
| 2.5e−3 | 3.06e−6 | 1.91 | 2.58e−6 | 1.97 |
| 1.25e−3 | 7.88e−7 | 1.96 | 6.67e−7 | 1.95 |
| 6.25e−4 | 1.99e−7 | 1.99 | 1.69e−7 | 1.98 |
| 3.125e−4 | 5.10e−8 | 1.96 | 4.25e−8 | 1.99 |
| 1.5625e−4 | 1.24e−8 | 2.04 | 1.06e−8 | 2.00 |
| 7.8125e−5 | 3.12e−9 | 1.99 | 2.66e−9 | 1.99 |

TABLE 4.1. The $L^2$ errors at $t = 0.1$ for the approximate phase variable $\phi$ of the presumed exact solution (4.2), computed by the schemes BDF2-L for the linear regularization model, and BDF2-W for the Willmore regularization model using various temporal resolutions.

| $\delta t$ | BDF2-L | Order | BDF2-W | Order |
|---|---|---|---|---|
| 2e−3 | 1.63e−5 | | 1.88e−5 | |
| 1e−3 | 5.29e−6 | 1.62 | 6.14e−6 | 1.60 |
| 5e−4 | 1.78e−6 | 1.57 | 2.08e−6 | 1.57 |
| 2.5e−4 | 5.78e−7 | 1.62 | 6.91e−7 | 1.58 |
| 1.25e−4 | 1.73e−7 | 1.74 | 1.85e−7 | 1.90 |
| 6.25e−5 | 4.11e−8 | 2.07 | 4.81e−8 | 1.94 |

TABLE 4.2. The $L^2$ numerical errors at $t = 0.1$ that are computed using schemes BDF2-L and BDF2-W using various temporal resolutions with the initial conditions of (4.3), for mesh refinement test in time.

4.1.1. *Presumed exact solution.* In the first example, we test the convergence rates of the proposed schemes for the isotropic model, i.e., $\gamma(\boldsymbol{n}) = 1$. We assume the following function

(4.2) $$\phi(x, y, t) = \sin(x)\cos(y)\cos(t)$$

to be the exact solution, and impose some suitable force fields such that the given solution can satisfy the linear regularization system (2.8)-(2.9).

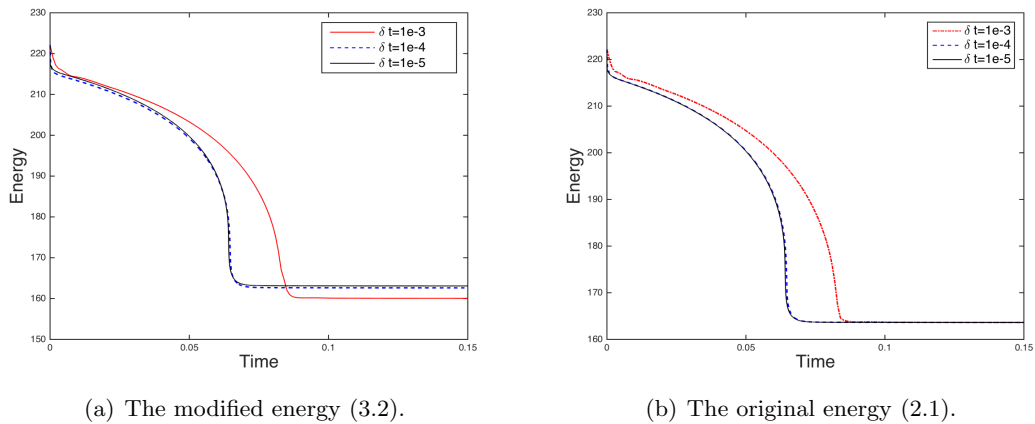

(a) The modified energy (3.2).  (b) The original energy (2.1).

FIGURE 4.1. Time evolution of the two free energy functionals for the isotropic model with linear regularization with $\beta = 6e-4$ and three different time steps $\delta t = $ 1e−3, 1e−4, 1e−5. (a) is the modified energy (3.2), and (b) is the modified energy (2.1). The energy curves show the decays for all time steps, which confirms that our algorithm is unconditionally stable.

We set $M(\phi) = \epsilon^2, \alpha = 0$ and all other parameters are from (4.1). In Table 4.1, we list the $L^2$ errors of the phase variable $\phi$ between the numerically simulated solution and the exact solution at $t = 0.1$ with different time step sizes. We observe that the schemes BDF2-L and BDF2-W achieve almost perfect second order accuracy in time.

4.1.2. *Mesh refinement in time.* We further perform the refinement test for the temporal convergence of the isotropic model. We set the initial conditions as follows,

$$(4.3) \qquad \phi(x, y, t = 0) = -\tanh\Big(\frac{\sqrt{(x-\pi)^2 + (y-\pi)^2} - 1.7}{2\epsilon}\Big).$$

Since the exact solutions are not known, we choose the solution obtained with the time step size $\delta t = 6.25e-5$ as the benchmark solution (approximately the exact solution) for computing errors. We present the $L^2$ errors of the phase variable between the numerical solution and the exact solution at $t = 0.1$ with different time step sizes in Table 4.2. We observe that the schemes asymptotically match the second order accuracy in time.

4.2. **Isotropic case with linear regularization.** In this subsection, we test the scheme BDF2-L for solving the isotropic case ($\gamma(\boldsymbol{n}) = 1$) with the linear regularization. We set the initial condition of two 2D circles as

$$(4.4) \qquad \phi(x, y, t = 0) = \sum_{i=1}^{2} -\tanh(\frac{\sqrt{(x-x_i)^2 + (y-y_i)^2} - r_i}{1.2\epsilon}) + 1,$$

where $(x_1, y_1, r_1) = (\pi - 0.7, \pi - 0.6, 1.5)$ and $(x_2, y_2, r_2) = (\pi + 1.65, \pi + 1.6, 0.7)$.

We emphasize that any time step size $\delta t$ is allowable for the computations from the stability concern since all developed schemes are unconditionally energy stable. But larger time step will definitely induce large numerical errors. Therefore, we need to discover the rough range of the allowable maximum time step size in order to obtain good accuracy and to consume as low computational cost as possible. This time step range could be estimated through the energy evolution curve plots, shown in Fig. 4.1, where we compare the time evolution of the free energy for three different time step sizes until the equilibrium using the scheme BDF2-L. We observe that all three energy curves decay monotonically for all time step sizes, which numerically confirms that our algorithms are unconditionally energy stable. For smaller time steps of $\delta t = 1\mathrm{e}{-5}$ and $1\mathrm{e}{-4}$, the two energy curves coincide very well. But for the larger time step of $\delta t = 1\mathrm{e}{-3}$, the energy curve deviates viewable away from others. This means the adopted time step size should not be larger than $1\mathrm{e}{-3}$, in order to get reasonably good accuracy.

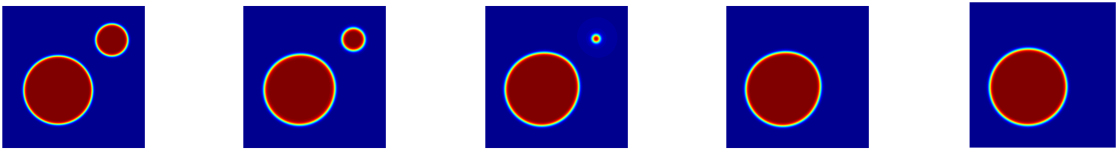

FIGURE 4.2. The 2D dynamical evolution of the phase variable $\phi$ with the initial condition (4.4) and time step $\delta t = 1\mathrm{e}{-4}$ for the isotropic model with the linear regularization with $\beta = 6\mathrm{e}{-4}$. Snapshots of the numerical approximation are taken at $t = 0$, $4\mathrm{e}{-2}$, $6.36\mathrm{e}{-2}$, $7\mathrm{e}{-3}$ and $2\mathrm{e}{-2}$.

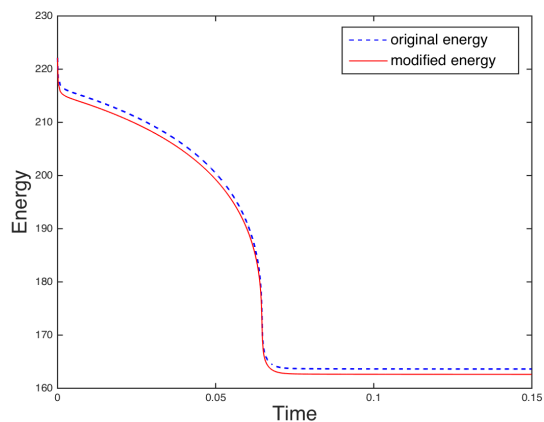

FIGURE 4.3. Time evolution of the two free energy functionals, the modified energy (3.2) and the original energy (2.1), for the isotropic model with the linear regularization by using $\beta = 6\mathrm{e}{-4}$ and initial condition (4.4).

In Fig. 4.2, we show the evolutions of the phase field variable $\phi$ at various time by using the time step $\delta t = 1\mathrm{e}{-4}$. we observe the coarsening effect that the small circle is absorbed into the big circle, and the total absorption happens around $t = 6.36\mathrm{e}{-2}$. We then present the evolution of the two free energy functionals, the modified energy (3.2) and original energy (2.1), in Fig. 4.3. It is remarkable that these two energies differ indeed, but both of them decay to the equilibrium, monotonically and simultaneously. At around $t = 6.36\mathrm{e}{-2}$, the energies undergo a rapid decrease when the totally absorption happens. Soon after that, the system achieves the equilibrium of circular shape immediately.

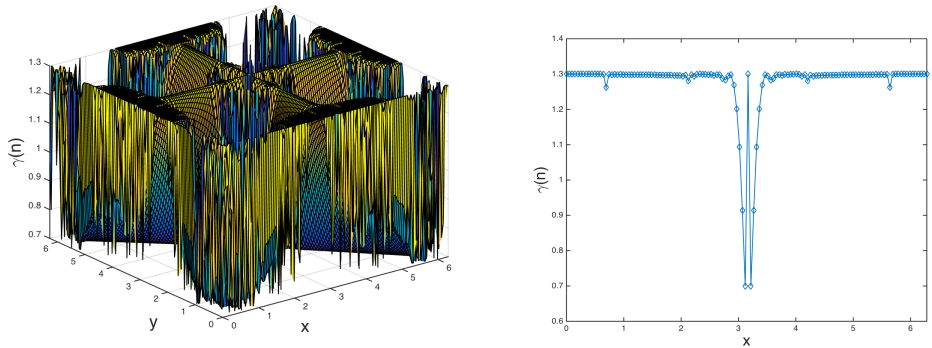

FIGURE 4.4. The profile of $\gamma(\boldsymbol{n}^0)$ with $\alpha = 0.3$ by using the initial condition (4.3). The left subfigure is the 2D surface plots of $\gamma(\boldsymbol{n}^0)$, and the right subfigure is the 1D cross-section of $\gamma(\boldsymbol{n}_0)|_{(\cdot, y=\pi)}$.

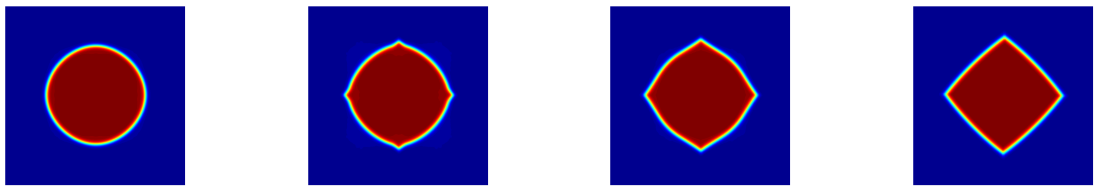

FIGURE 4.5. The 2D dynamical evolution of the phase variable $\phi$ for the anisotropic linear regularization model by using the initial condition (4.3), the time step $\delta t = 1\mathrm{e}{-4}$, $\alpha = 0.3$, and $\beta = 1\mathrm{e}{-4}$. Snapshots of the numerical approximation are taken at $t = 5\mathrm{e}{-4}$, $1.5\mathrm{e}{-3}$, $2.5\mathrm{e}{-3}$ and $2\mathrm{e}{-2}$.

4.3. **Anisotropic linear regularization model.** In this subsection, we consider the anisotropic system with linear regularization. We set $M(\phi) = 1$ and all other parameters are from (4.1), if not explicitly specified.

4.3.1. *Evolution of a 2D circle.* We perform the simulation for the evolution of a circle in 2D by using the initial condition $\phi^0$ specified in (4.3). In Fig. 4.4, we present the 2D profile of $\gamma(\boldsymbol{n}_0)$ and 1D cross-section of $\gamma(\boldsymbol{n}_0)|_{(\cdot, y=\pi)}$, where a high oscillation profile is found almost everywhere.

In Fig. 4.5, we show the dynamics how a circular shape interface with full orientations evolves to an anisotropic pyramid with missing orientations at four corners. Snapshots of the phase field variable $\phi$ are taken at $t = 5\mathrm{e}{-4}, 1.5\mathrm{e}{-3}, 2.5\mathrm{e}{-3}$ and $2\mathrm{e}{-2}$. In Fig. 4.6, we show the evolution of the two free energy functionals until the steady state for the modified energy (3.2) and original energy (2.1).

In Fig. 4.7, we show the evolution of the original free energy functional (2.1) for four combinations of stabilizers: (i) $S_1 = S_2 = S_3 = 0$; (ii) $S_1 = S_2 = 0, S_3 = 0.001$; (iii) $S_1 = S_2 = 2, S_3 = 0$; and

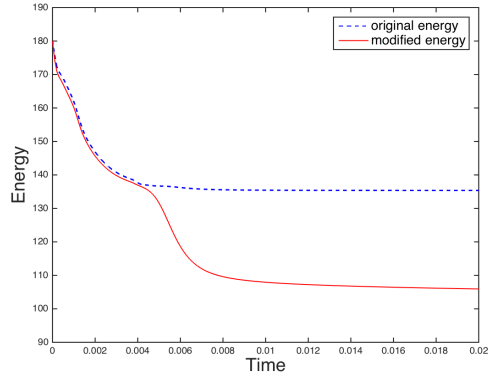

FIGURE 4.6. Time evolution of the two free energy functionals, the modified energy (3.2) and the original energy (2.1), for the anisotropic model with linear regularization with $\alpha = 0.3, \beta = 6e{-}4$ and initial condition (4.3).

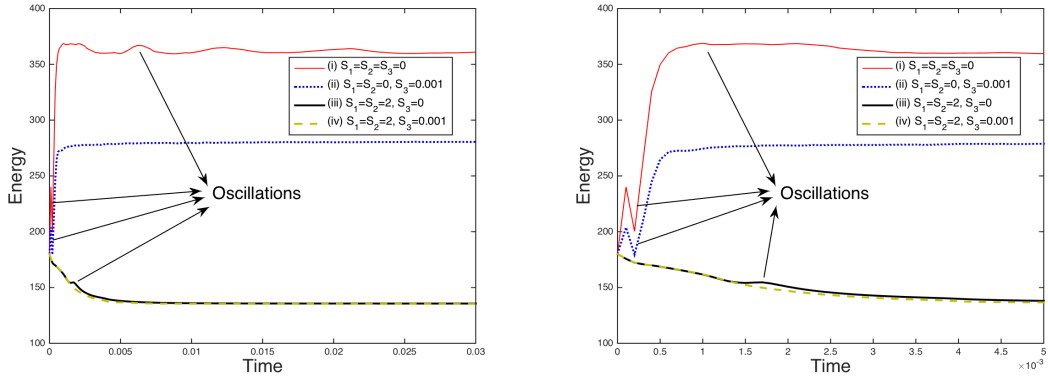

(a) Energy evolutions with four combinations of stabilizers.

(b) A close-up view.

FIGURE 4.7. Time evolution of the original free energy functional (2.1) of four combinations of linear stabilizers, for the anisotropic model with the linear regularization. The initial condition is (4.3), $\alpha = 0.3$ and $\beta = 6e{-}4$. The left subfigure (a) is the energy profile for $t \in [0, 3e{-}2]$, and the right subfigure (b) is a close-up view for $t \in [0, 3e{-}3]$.

(iv) $S_1 = S_2 = 2$, $S_3 = 0.001$. For (i) and (ii), the energies present some non-physical oscillations and further increase with time; for (iii), the energy mostly decays but still presents a very slight oscillation around $t = 0.002$; and for (iv), all oscillations vanish and energy decays monotonically, that means the combination (iv) which is adopted for simulations, is a combination of stabilizers that can suppress high-frequency oscillations efficiently.

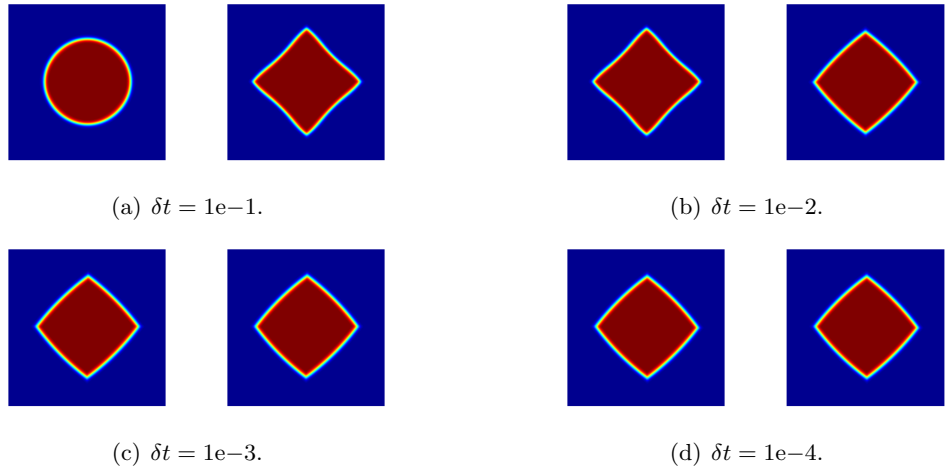

(a) $\delta t = 1\text{e}{-}1$.  (b) $\delta t = 1\text{e}{-}2$.

(c) $\delta t = 1\text{e}{-}3$.  (d) $\delta t = 1\text{e}{-}4$.

FIGURE 4.8. Snapshots of the phase variable $\phi$ for four different time steps, where (a) $\delta t = 1\text{e}{-}1$; (b) $\delta t = 1\text{e}{-}2$; (c) $\delta t = 1\text{e}{-}3$; and (d) $\delta t = 1\text{e}{-}4$. The initial condition is (4.3), $\alpha = 0.3$, and $\beta = 6\text{e}{-}4$. For each panel, the left subfigure is the profile of $\phi$ at $t = 0.2$, and the right subfigure is that of $t = 2$.

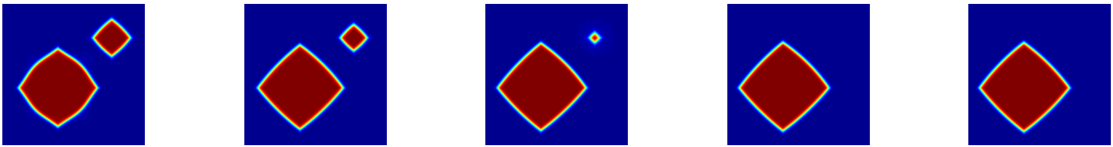

FIGURE 4.9. The 2D dynamical evolution of the phase variable $\phi$ for the anisotropic model with the linear regularization. The initial condition is (4.4), $\delta t = 1\text{e}{-}4$, $\alpha = 0.3$, and $\beta = 6\text{e}{-}4$. Snapshots of the numerical approximation are taken at $t = 2\text{e}{-}3, 6\text{e}{-}2, 9.2\text{e}{-}2, 1\text{e}{-}1$, and $2\text{e}{-}1$.

We further testify if the scheme BDF2-L can allow large time steps. In Fig. 4.8, we show the snapshots of the phase field variable $\phi$ at $t = 0.2$ and 2 for four different time steps $\delta t = 1\text{e}{-}1, 1\text{e}{-}2, 1\text{e}{-}3$ and $1\text{e}{-}4$. We observe that, for the two smaller time steps $\delta t = 1\text{e}{-}3$ and $1\text{e}{-}4$, the profiles coincide very well. But for the larger time step of $\delta t = 1\text{e}{-}1$ and $1\text{e}{-}2$, the results deviate viewable away from others. This means the so-called "unconditional energy stability" here means the schemes have no constraints for the time step only from the stability concern. It does not mean that any arbitrarily large time step can be chosen for computations since larger time step size will definitely induce larger numerical errors in practice.

4.3.2. *Evolution of two 2D circles.* In this example, we use the initial condition given in (4.4) of two circles to see how the combined effects of anisotropy and coarsening execute. In Fig. 4.9, snapshots of the profiles of the phase field variable $\phi$ are taken at $t = 2\text{e}{-}3, 6\text{e}{-}2, 9.2\text{e}{-}2, 1\text{e}{-}1$ and $2\text{e}{-}1$. We observe that the two circles first evolve to anisotropic shapes with missing orientation at the four corners, then the anisotropic system coarsens and the small shape disappears. In

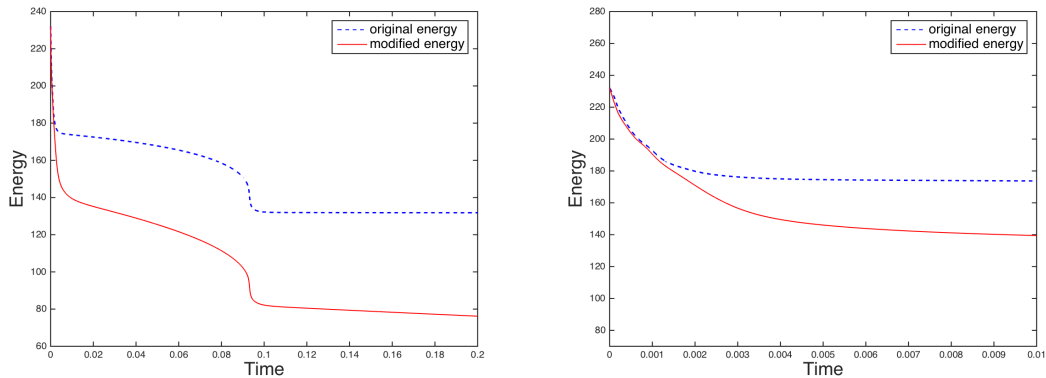

FIGURE 4.10. Time evolution of the two free energy functionals, the modified energy (3.2) and the original energy (2.1), for the anisotropic model with the linear regularization, by using the initial condition (4.4), $\delta t = 1\text{e}-4$, $\alpha = 0.3$, and $\beta = 6\text{e}-4$. The left subfigure is the energy profile until the equilibrium, and the right subfigure is a close-up showing where the energy decreases fast.

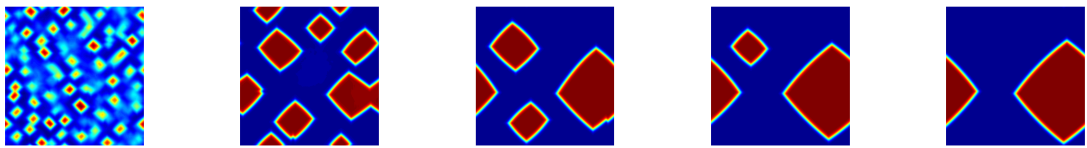

FIGURE 4.11. The 2D dynamical evolution of the phase variable $\phi$ of the spinodal decomposition example for the anisotropic linear regularization model, by using the initial condition (4.5), $\delta t = 1\text{e}-4$, $\alpha = 0.3$, and $\beta = 6\text{e}-4$. Snapshots are taken at $t = 1.5\text{e}-3, 1\text{e}-2, 2.5\text{e}-1, 5\text{e}-1, 1$.

Fig. 4.10, we present the evolution of the free energy functional until the steady state for the modified energy (3.2) and original energy (2.1). The energies undergo a rapid decrease when the totally absorption happens at around $t = 9.2\text{e}-2$. Then the system achieves the equilibrium of circular shape immediately after that.

4.3.3. *Spinodal decomposition in 2D.* In this example, we study the phase separation dynamics that is called spinodal decomposition for the anisotropic model with the linear regularization using the scheme BDF2-L. The process of the phase separation can be studied by considering a homogeneous binary mixture, which is quenched into the unstable part of its miscibility gap. When the spinodal decomposition takes place, the spontaneous growth of the concentration fluctuations are manifested that leads the system from the homogeneous to the two-phase state. The initial condition is taken as the randomly perturbed concentration field as follows,

(4.5) $$\phi(x, y, t = 0) = -0.4 + 0.001\text{rand}(x, y).$$

In Fig. 4.11, we perform the simulations by using the time step $\delta t = 1\text{e}-4$. We observe the combined effects of anisotropy and coarsening when time evolves. The final equilibrium solution is obtained after $t = 1$, where equilibrium shape becomes a pyramid (independent to each other notably) due to the strong anisotropy. In Fig. 4.12, we plot the evolution of the modified energy (3.2) and the original free energy (2.1), which show the decays with time that confirms that our algorithms are unconditionally stable.

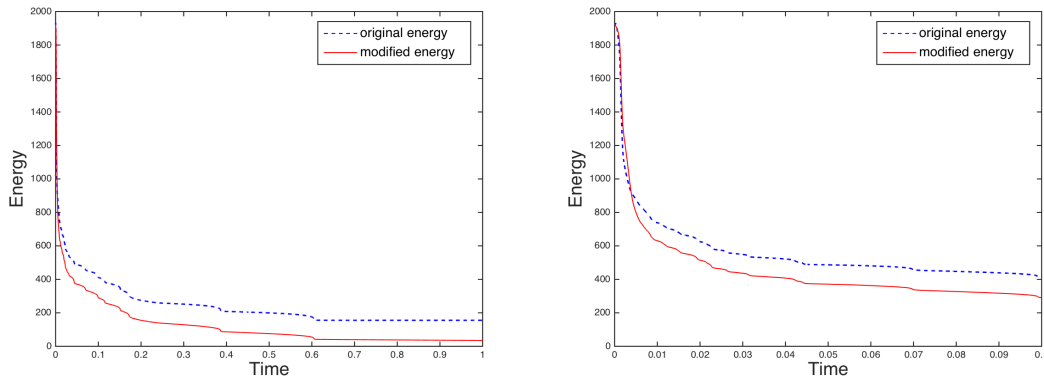

FIGURE 4.12. Time evolution of the two free energy functionals, the modified energy (3.2) and the original energy (2.1), of the spinodal decomposition example for the anisotropic linear regularization model by using the initial condition (4.5), $\delta t = 1\text{e}{-4}$, $\alpha = 0.3$, and $\beta = 6\text{e}{-4}$. The left subfigure is the energy profile until the equilibrium, and the right subfigure is a close-up view showing where the energy decreases fast.

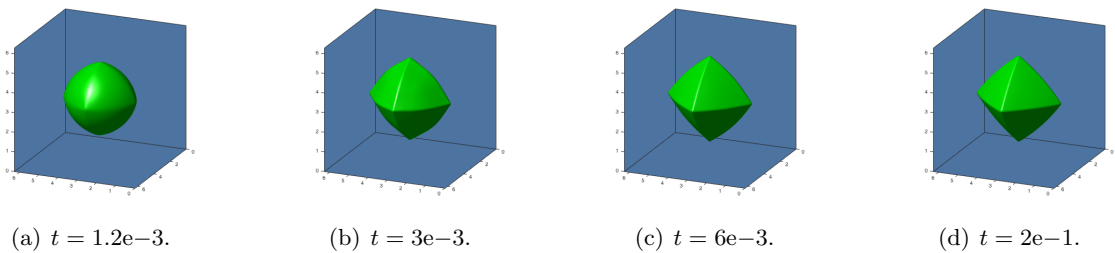

(a) $t = 1.2\text{e}{-3}$.  (b) $t = 3\text{e}{-3}$.  (c) $t = 6\text{e}{-3}$.  (d) $t = 2\text{e}{-1}$.

FIGURE 4.13. The dynamical evolution of a 3D sphere for the anisotropic linear regularization model, by using the initial condition (4.6), $\delta t = 1\text{e}{-4}$, $\alpha = 0.3$, and $\beta = 6\text{e}{-4}$. Snapshots of the isosurfaces of the phase field variable $\{\phi = 0\}$ are taken at $t = 1.2\text{e}{-3}, 3\text{e}{-3}, 6\text{e}{-3}$ and $2\text{e}{-1}$.

4.3.4. *Evolution of a 3D sphere.* We next investigate the 3D simulation of a sphere by using the following initial condition

(4.6) $$\phi(x, y, z, t = 0) = -\tanh\Big(\frac{\sqrt{(x-\pi)^2 + (y-\pi)^2 + (z-\pi)^2} - 1.7}{2\epsilon}\Big).$$

The evolution of the spherical shape towards its equilibrium is shown in Fig. 4.13, where we observe that the 3D sphere evolves to an anisotropic pyramid with missing orientations at six corners.

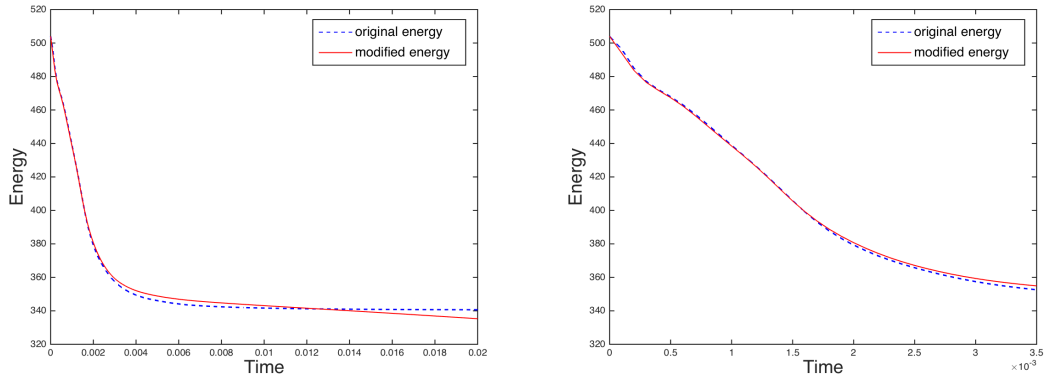

FIGURE 4.14. Time evolution of the two free energy functional for the 3D sphere example, the modified energy (3.2) and the original energy (2.1), for the anisotropic model with linear regularization with $\alpha = 0.3$, $\beta = 6e-4$ and initial condition (4.6). The left subfigure is the energy profile until the equilibrium, and the right subfigure is a close-up showing where the energy decreases fast.

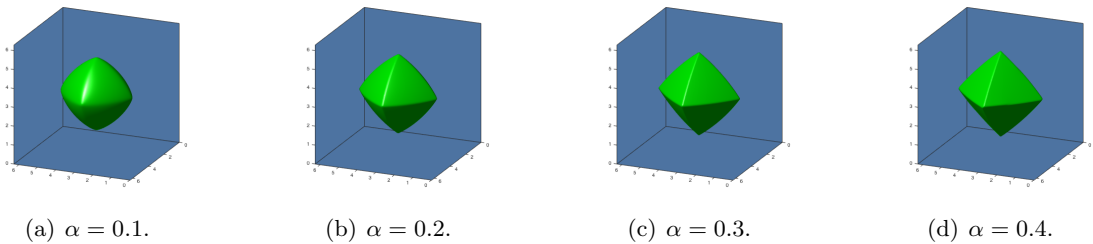

(a) $\alpha = 0.1$.    (b) $\alpha = 0.2$.    (c) $\alpha = 0.3$.    (d) $\alpha = 0.4$.

FIGURE 4.15. Snapshots of the equilibrium solution of the anisotropic linear regularization model at $t = 0.03$ for four different anisotropic parameters $\alpha$ (from left to right): 0.1, 0.2, 0.3, and 0.4, respectively, by using the initial condition (4.6), $\beta = 6e-4$, $\delta t = 1e-4$.

Snapshots of the profiles of the isosurfaces of $\{\phi = 0\}$ are taken at $t = 1.2e-3$, $3e-3$, $6e-3$ and $2e-1$. In Fig. 4.14, we present the evolution of the free energy functional until the steady state for the modified energy (3.2) and original energy (2.1).

We further investigate the effect of the strength of the anisotropy parameter $\alpha$ on crystal shapes in 3D. We choose four different values of $\alpha = 0.1, 0.2, 0.3, 0.4$ and list the steady state solutions for these four cases in Fig. 4.15. The corresponding 2D cross-sections along $x = \pi$ are shown in Fig. 4.16, respectively. When $\alpha$ is relatively small, e.g., $\alpha = 0.1$, missing orientations are not seen clearly, see Fig. 4.15 (a). When $\alpha$ increases (e.g., 0.2, 0.3, 0.4), equilibrium shapes tend to become pyramids due to the strong anisotropy, see Fig. 4.15 (b, c, d). In Fig. 4.17, we show the evolution of the original free energy (2.1), we can observe that the energy decreases faster with larger $\alpha$. The results are in very good agreement with simulations shown in [2, 21].

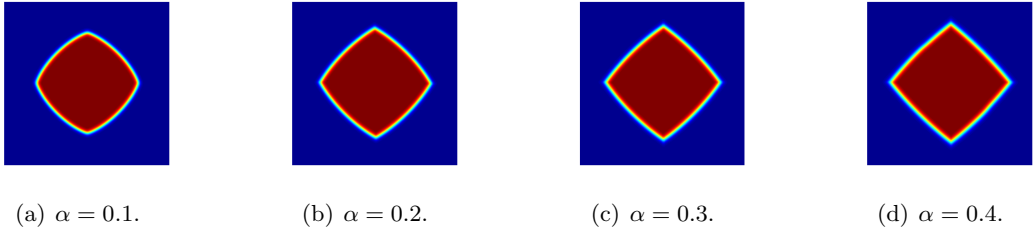

(a) $\alpha = 0.1$.    (b) $\alpha = 0.2$.    (c) $\alpha = 0.3$.    (d) $\alpha = 0.4$.

FIGURE 4.16. 2D cross-section of $\phi(\pi, \cdot, \cdot)$ of the equilibrium solution from Fig. 4.15.

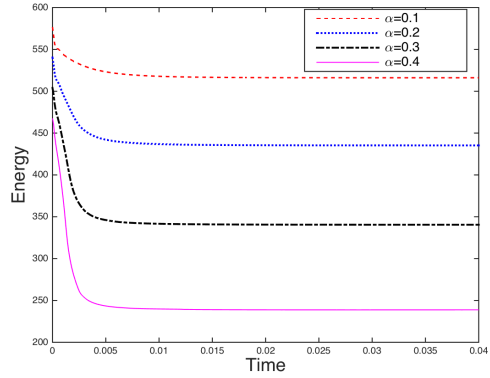

FIGURE 4.17. Time evolutions of the free energy functional (2.1) for the anisotropic model with the linear regularization, by using four different values of $\alpha$ shown in Fig. 4.15 and Fig. 4.16.

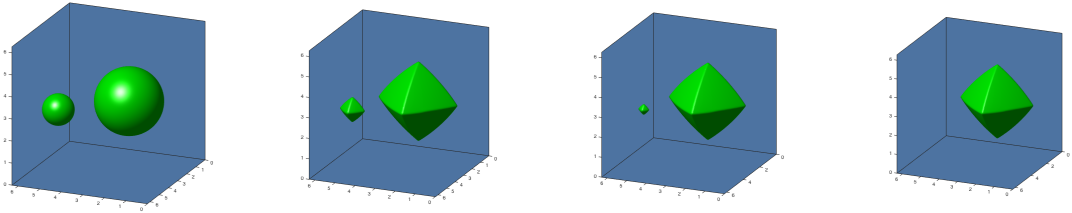

FIGURE 4.18. The dynamical evolution of two 3D spheres for the anisotropic model with the linear regularization, by using the initial condition (4.7), $\alpha = 0.3$, $\beta = 6e-4$, and $\delta t = 1e-4$. Snapshots of the isosurfaces of the phase field variable $\phi$ are taken at $t = 0$, $1.32e-2$, $1.98e-2$ and $4e-2$.

4.3.5. *Evolution of two 3D spheres.* In this example, we implement the 3D simulation of two spheres for the anisotropic linear regularization model by using the following initial condition

$$(4.7) \qquad \phi(x, y, t = 0) = \sum_{i=1}^{2} -\tanh\left(\frac{\sqrt{(x-x_i)^2 + (y-y_i)^2 + (z-z_i)} - r_i}{1.2\epsilon}\right) + 1,$$

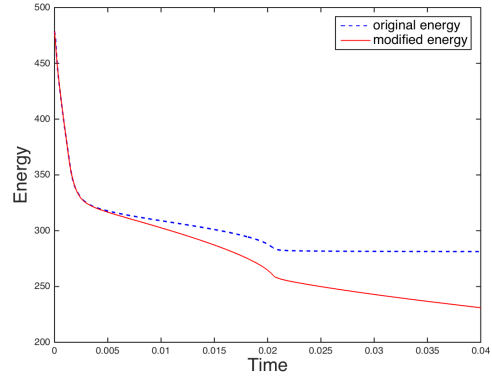

FIGURE 4.19. Time evolution of the two free energy functionals of two 3D spheres example, the modified energy (3.2) and the original energy (2.1), for the anisotropic model with linear regularization by using the initial condition (4.7), $\alpha = 0.3$, $\beta = 6e-4$, and $\delta t = 1e-4$.

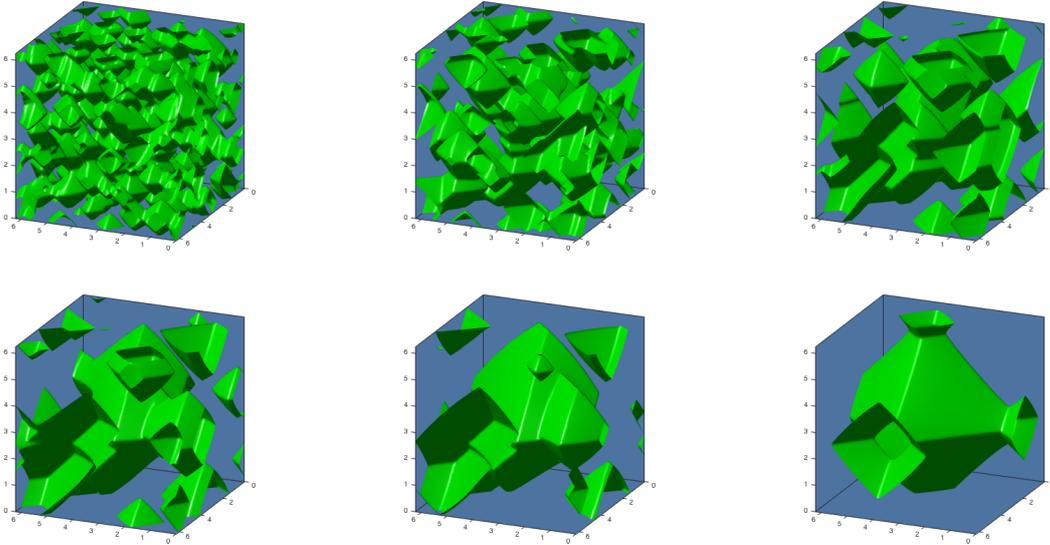

FIGURE 4.20. The dynamical evolution of the phase variable $\phi$ of the 3D spinodal decompositions for the anisotropic linear regularization model, by using the initial condition (4.8), $\delta t = 1e-4$, $\alpha = 0.3$, and $\beta = 6e-4$. Snapshots of the isosurfaces of $\{\phi = 0\}$ are taken at $t = 1e-3, 3e-2, 1e-1, 2e-1, 3e-1$ and $1$.

where $(x_1, y_1, z_1, r_1) = (\pi - 0.7, \pi - 0.6, \pi, 1.5)$ and $(x_2, y_2, z_2, r_2) = (\pi + 1.65, \pi + 1.6, \pi, 0.7)$.

The evolution of the two 3D spheres towards its equilibrium is shown in Fig. 4.18. Snapshots of the profiles of the isosurfaces of $\{\phi = 0\}$ are taken at $t = 0, 1.32e-2, 1.98e-2$ and $4e-2$. The two spheres first evolve to anisotropic pyramids with missing orientation at the six corners, then the

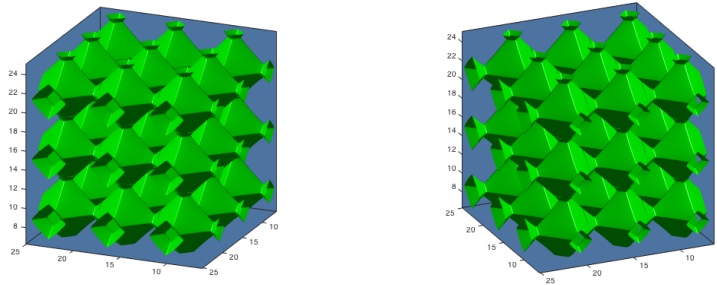

FIGURE 4.21. The isosurfaces of the steady state solution at $t = 1$ of the 3D spinodal decompositions for the anisotropic linear regularization model for 4 periods, i.e., $[0, 8\pi]^3$. The two subfigures are from different angles of view.

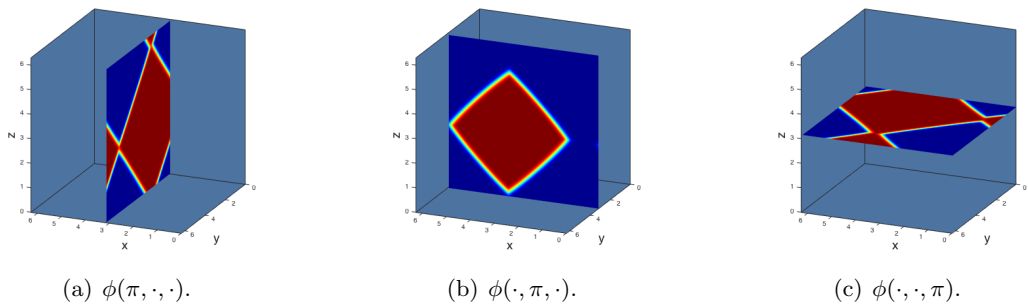

(a) $\phi(\pi, \cdot, \cdot)$.  (b) $\phi(\cdot, \pi, \cdot)$.  (c) $\phi(\cdot, \cdot, \pi)$.

FIGURE 4.22. 2D cross-section for the steady state solution at $t = 1$ from Fig. 4.20. From left to right: $\phi(\pi, \cdot, \cdot)$, $\phi(\cdot, \pi, \cdot)$ and $\phi(\cdot, \cdot, \pi)$, respectively.

anisotropic system coarsens and the small pyramid is totally absorbed into the larger one. The 3D dynamics are consistent with the 2D example (Fig. 4.9). In Fig. 4.19, we present the evolution of the free energy functional until the steady state for the modified energy (3.2) and original energy (2.1).

4.3.6. *Spinodal decomposition in 3D.* In this example, we study the 3D spinodal decomposition for the anisotropic model with linear regularization using the scheme BDF2-L. To be consistent with the 2D example, we use the initial condition as follows,

$$(4.8) \qquad \phi(x, y, z, 0) = -0.4 + 0.001\text{rand}(x, y, z).$$

In Fig. 4.20, we present evolutions of the phase field variable $\phi$ by showing the snapshots of the isosurfaces of $\{\phi = 0\}$. We still take the time step $\delta t = 1\text{e}{-4}$ and the final equilibrium solution is obtained around $t = 1$. The equilibrium solution presents pyramid-like shapes due to the strong anisotropy. To see the equilibrium shape more clearly, since the computed domain is periodic, we then plot the isosurfaces for 4 periods, i.e., $[0, 8\pi]^3$ in Fig. 4.21 using two different angles of view. We observe that all pyramids are actually connected together, that is quite different from the 2D example where the pyramids are located independently (cf. the last subfigure in Fig. 4.11). But in Fig. 4.22, when we observe the 2D cross-sections of the equilibrium solution, we can find that the cross-section $\phi(\cdot, \pi, \cdot)$ is the same as the 2D results. These findings actually imply that the 3D results are essentially different from 2D, and 3D computations can not be simply predicted by 2D as well. In Fig. 4.23, we further show the evolution of the modified energy (3.2) and the original free energy (2.1).

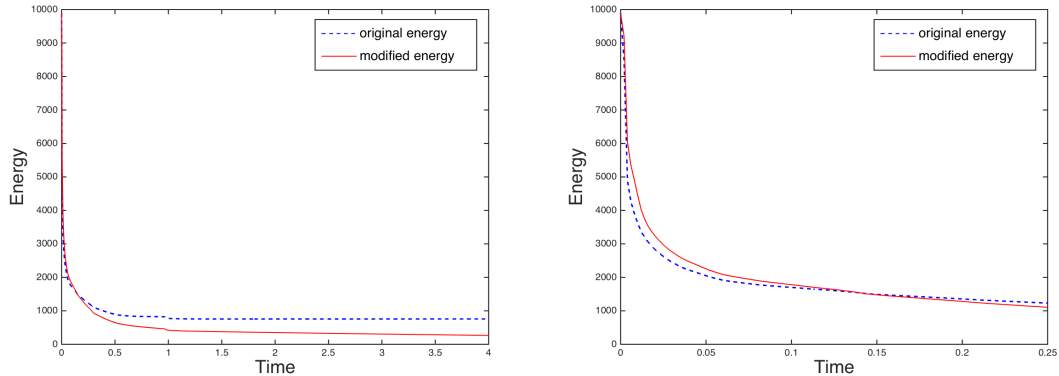

FIGURE 4.23. Time evolution of the two free energy functionals for two 3D spheres example, the modified energy (3.2) and the original energy (2.1), for the anisotropic model with linear regularization with the initial condition (4.8), $\alpha = 0.3$, $\beta = 6\mathrm{e}{-4}$. The left subfigure is the energy profile until the equilibrium, and the right subfigure is a close-up showing where the energy decreases fast.

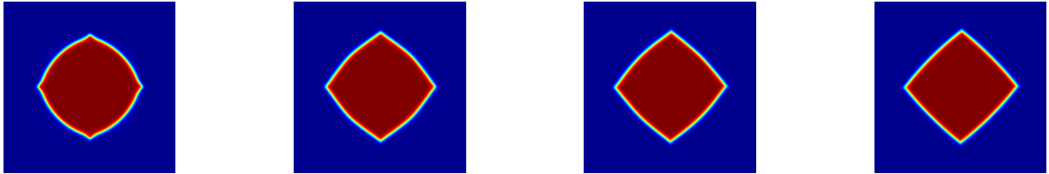

FIGURE 4.24. The 2D dynamical evolution of the phase variable $\phi$ for the anisotropic Willmore regularization model by using the initial condition (4.3), $\delta t = 1\mathrm{e}{-4}$, $\alpha = 0.3$, $\beta = 6\mathrm{e}{-4}$. Snapshots are taken at $t = 2\mathrm{e}{-3}$, $4\mathrm{e}{-3}$, $6\mathrm{e}{-3}$ and $3\mathrm{e}{-2}$.

**4.4. Anisotropic Willmore regularization model.** We consider the anisotropic system with the Willmore regularization in this subsection. To be consistent with the linear regularization model, we still let $M(\phi) = 1$ and use all other parameters are from (4.1), if not explicitly specified.

We implement the simulation for the evolution of a 2D circle by using the scheme BDF2-W and the initial condition (4.3). In Fig. 4.24, we show how an the circular interface with full orientations evolves to an anisotropic one with missing orientations at four corners. Snapshots of the phase field variable $\phi$ are taken at $t = 2\mathrm{e}{-3}$, $4\mathrm{e}{-3}$, $6\mathrm{e}{-3}$ and $3\mathrm{e}{-2}$. In Fig. 4.25, we plot the free energy functionals for the modified energy (3.2) and original energy (2.1) until the equalibrium.

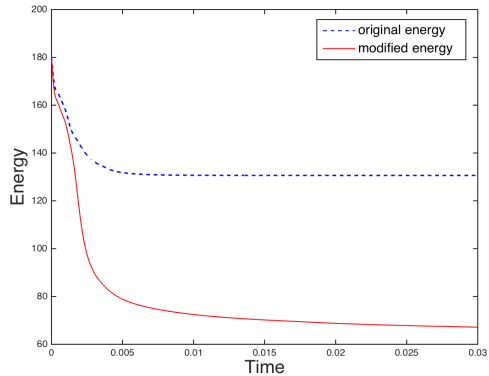

FIGURE 4.25. Time evolution of the two free energy functionals for the 2D circle example, the modified energy (3.2) and the original energy (2.1), for the anisotropic model with the Willmore regularization by using the initial condition (4.3), $\alpha = 0.3$, $\beta = 6e-4$.

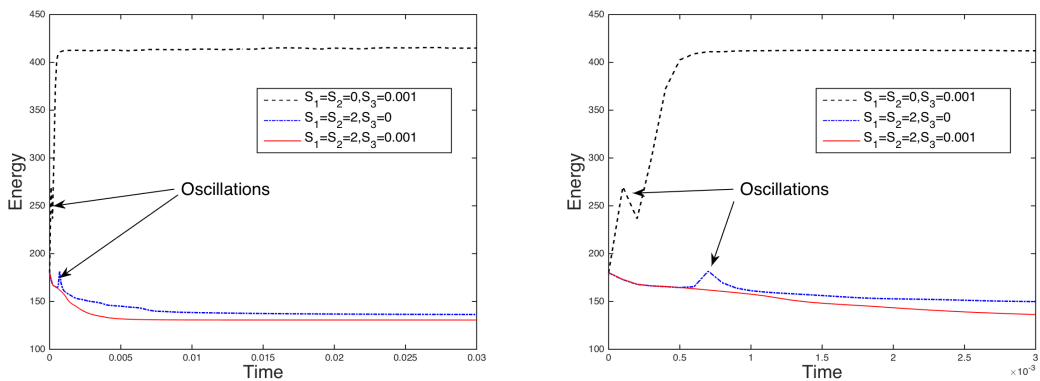

FIGURE 4.26. Time evolution of the original energy (2.1) when using three combinations of linear stabilizers, for solving the anisotropic model with the Willmore regularization. The results are computed by using the initial condition (4.3), $\alpha = 0.3$, and $\beta = 6e-4$. The left subfigure is the energy profile for $t \in [0, 3e-2]$, and the right subfigure is a close-up view for $t \in [0, 3e-3]$.

In Fig. 4.26, we present the evolution of the original free energy functional (2.1) for three combinations of stabilizers: (i) $S_1 = S_2 = 0$, $S_3 = 0.001$; (ii) $S_1 = S_2 = 2$, $S_3 = 0$; and (iii) $S_1 = S_2 = 2$, $S_3 = 0.001$. As the linear regularization case, combinations (i) and (ii) cause non-physical oscillations, and combination (iii) can suppress high-frequency oscillations efficiently.

We further testify if our developed scheme BDF2-W can allow large time steps. In Fig. 4.27, we show the snapshots of the phase field variable $\phi$ at $t = 0.2$ and 2 for four different time steps $\delta t = 1e-1, 1e-2, 1e-3$, and $1e-4$. We observe that, for the three smaller time steps (e.g. $1e-2$, $1e-3$ and $1e-4$), the profiles coincide very well. But for the largest time step of $1e-1$, the results deviate viewable away from others. Therefore, we choose $\delta t = 1e-4$, which is made by the accuracy requirement instead of the stability requirement.

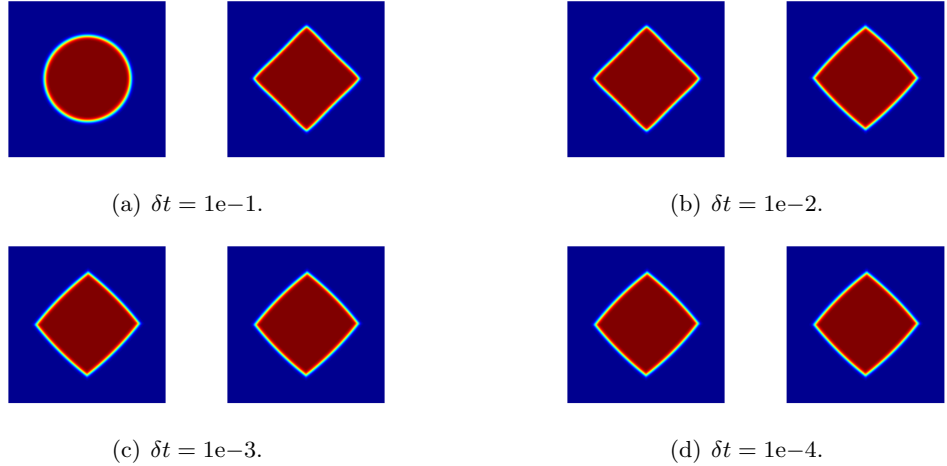

(a) $\delta t = 1e-1$.   (b) $\delta t = 1e-2$.

(c) $\delta t = 1e-3$.   (d) $\delta t = 1e-4$.

FIGURE 4.27. Snapshots of the phase variable $\phi$ with four different time steps for solving the anisotropic model with the Willmore regularization, where (a) $\delta t = 1e-1$; (b) $\delta t = 1e-2$; (c) $\delta t = 1e-3$; and (d) $\delta t = 1e-4$. The initial condition is (4.3), $\alpha = 0.3$, and $\beta = 6e-4$. For each panel, the left subfigure is the profile of $\phi$ at $t = 0.2$, and the right subfigure is that of $t = 2$.

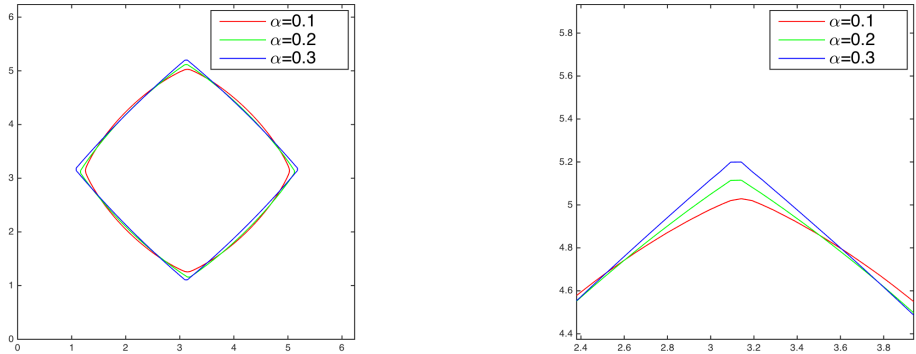

FIGURE 4.28. The contours of $\{\phi = 0\}$ of the equilibrium solutions for the anisotropic model with the Willmore regularization by using the initial condition (4.3), $\beta = 6e-4$ and three different values $\alpha = 0.1, 0.2, 0.3$.

Finally, we investigate the impacts of the strength of the anisotropy parameter $\alpha$ and the Willmore regularization parameter $\beta$ on the equilibrium four-fold shapes. We choose three different values of $\alpha = 0.1, 0.2, 0.3$ and fix $\beta = 6e-4$, the equilibrium shapes are compared in Fig. 4.28. It can be observed that the increase of $\alpha$ leads to a pyramid with shaper corners because of the strong anisotropy. The evolutions of the original total free energy (2.1) are shown in Fig. 4.29, where the energy decreases faster with larger $\alpha$. The effects of the Willmore regularization parameter $\beta$ on the equilibrium shapes are shown in Fig. 4.30, in which, we choose $\beta = 8e-4$, $1e-3$, $3e-3$ and fix $\alpha = 0.3$. It can be seen, as $\beta$ decreases, the corners become more shaper in the equilibrium morphologies. The corresponding energy evolutions are plotted in Fig. 4.31, where one can observe that energy decreases faster with smaller $\beta$. These numerical results are in very good agreement with the results shown in [2, 17, 21].

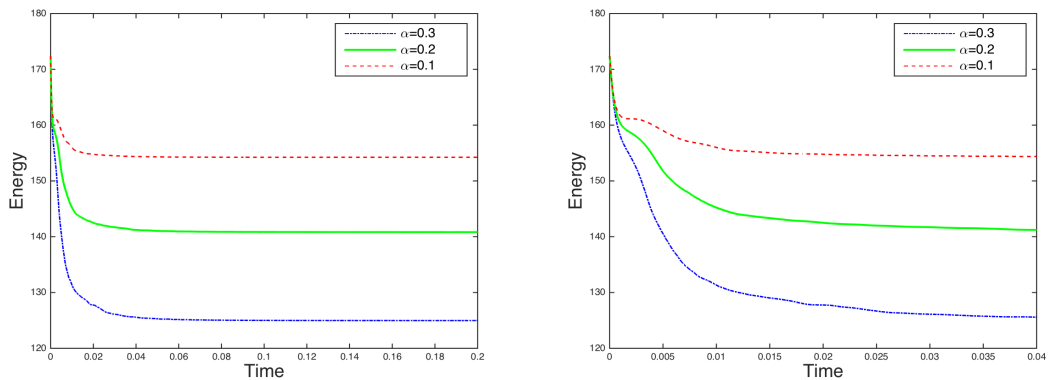

Figure 4.29. Time evolution of the original free energy functional (2.1) for the 2D circle example of the anisotropic model with the Willmore regularization, by using the initial condition (4.3), $\beta = 6e-4$ and three different values $\alpha = 0.1, 0.2, 0.3$.

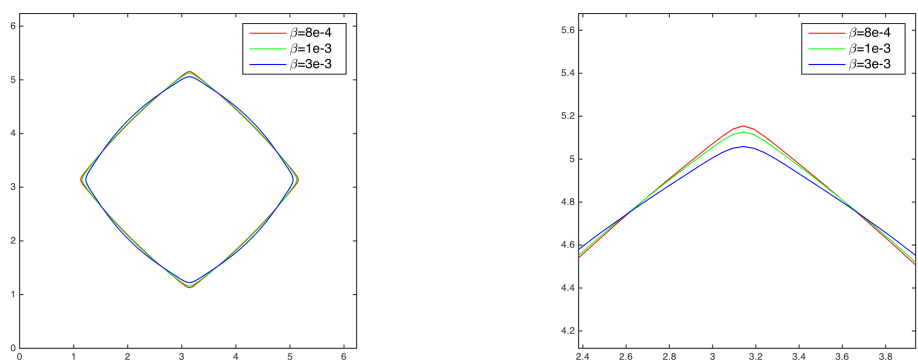

Figure 4.30. The contours of $\{\phi = 0\}$ of the equilibrium solutions for the anisotropic model with the Willmore regularization by using the initial condition (4.3), $\alpha = 0.3$ and three different values $\beta = 8e-4$, $1e-3$, $3e-3$.

## 5. Concluding Remarks

In this paper, we have developed two efficient, semi-discrete in time, second order linear schemes for solving the anisotropic Cahn-Hilliard phase field model by combining the linear stabilization approach and the recently developed SAV approach. The novelty of the proposed schemes is that all nonlinear terms can be treated semi-explicitly, and one only needs to solve three decoupled, sixth order equations with constant coefficients at each time step. We add three linear stabilization terms, which are shown to be crucial to remove the oscillations caused by the anisotropy coefficients. Compared to the existed schemes for the anisotropic model, our proposed schemes that conquer the inconvenience from nonlinearities by linearizing the nonlinear terms in the new way, are provably unconditionally energy stable, and thus allow for large time steps. We further numerically verify the accuracy in time and present various 2D and 3D numerical results for some benchmark numerical simulations.

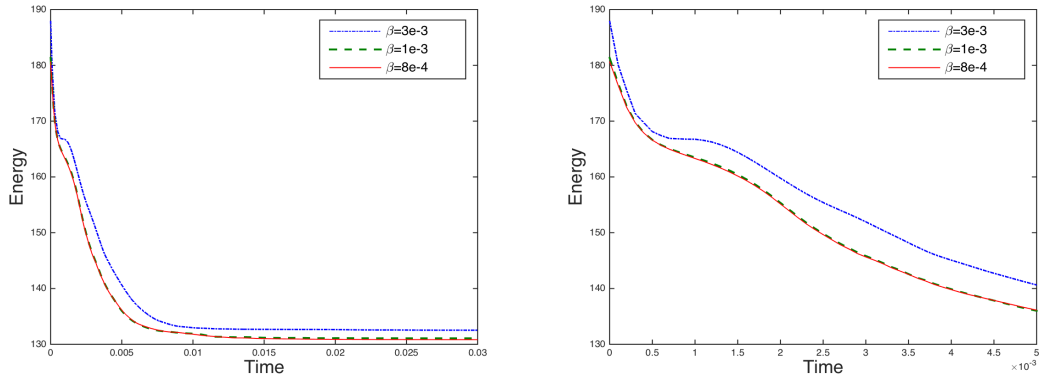

FIGURE 4.31. Time evolution of the original free energy functional (2.1) for the 2D circle example of the anisotropic model with the Willmore regularization, by using the initial condition (4.3), $\alpha = 0.3$ and three different values of $\beta = 8\mathrm{e}{-4}$, $1\mathrm{e}{-3}$, $3\mathrm{e}{-3}$.

Finally, the method is general enough to be extended to develop linear schemes for a large class of gradient flow problems with complex nonlinearities in the free energy density. Although we consider only time discrete schemes in this study, the results can be carried over to any consistent finite-dimensional Galerkin approximations since the proofs are all based on a variational formulation with all test functions in the same space as the space of the trial functions.